\documentclass[12pt]{article}
\usepackage{amsmath,amsthm,amsfonts,amscd,eucal}
\begin{document}
\baselineskip=20pt
\topmargin=2.0cm
\def\RR{{\mathbb R}}
\def\CC{{\mathbb C}}
\def\ZZ{{\mathbb Z}}
\def\sl2{{{\rm SL}(2,\RR)}}
\def\A{{\cal A}}
\def\B{{\cal B}}
\def\C{{\cal C}}
\def\D{{\cal D}}
\def\F{{\cal F}}
\def\H{{\cal H}}
\def\K{{\cal K}}
\def\k{{\rm K}}
\def\M{{\cal M}}
\def\N{{\cal N}}
\def\O{{\cal O}}
\def\P{{\cal P}}
\def\R{{\cal R}}
\def\V{{\cal V}}
\def\W{{\cal W}}
\def\G{{\bf G}}

\newtheorem{definition}{Definition}[section]
\newtheorem{theorem}{Theorem}[section]
\newtheorem{corollary}{Corollary}[section]
\newtheorem{proposition}{Proposition}[section]
\newtheorem{lemma}{Lemma}[section]

\title{Classification of Subsystems for Local Nets with
Trivial Superselection Structure\\
\bigskip
\small{
{\it Dedicated to S. Doplicher and J. E. Roberts on the
occasion of their 60th birthday}}
}

\author{Sebastiano Carpi$^*$, Roberto Conti$^{**}$\\
{}\\
$^*$ Dipartimento di Matematica \\
Universit\`a di Roma ``La Sapienza''\\
Piazzale A. Moro\\
I-00185 Roma, Italy\\
{}\\
$^{**}$ Dipartimento di Matematica \\
Universit\`a di Roma ``Tor Vergata''\\
Via della Ricerca Scientifica \\
I-00133 Roma, Italy}
\date{\today}
\maketitle
\markboth{S. Carpi, R. Conti}{On subsystems}
\renewcommand{\sectionmark}[1]{}
\begin{abstract} 
Let $\F$ be a local net of von Neumann algebras 
in four spacetime dimensions
satisfying certain natural structural assumptions.
We prove that if $\F$ has trivial superselection structure
then every covariant, Haag-dual subsystem $\B$
is of the form $\F_1^G \otimes I$ 
for a suitable decomposition $\F=\F_1 \otimes \F_2$ 
and a compact group action.
Then we discuss some application of our result, 
including free field models and certain 
theories with at most countably many sectors.
\end{abstract}

\vfill
\thanks{E-mail: carpi@mat.uniroma1.it, conti@mat.uniroma2.it.\par
$^{**}$ Supported by EU.}

\newpage

\section{Introduction}
In the algebraic approach to QFT \cite{Ha} the main objects under
investigation are 
(isotonous)
nets of von Neumann algebras over bounded regions in the Minkowski 
spacetime, satisfying pertinent additional requirements.
Any such correspondence is usually denoted by $\O \to \F(\O)$.

Internal symmetries of a net $\F$ can be defined as those automorphisms 
of the $C^*$-inductive limit $(\cup_{\O\in\K}\F(\O))^{-\|\cdot\|}$
(the {\it quasi-local} $C^*$-algebra;
it is customary to denote it in the same way as the net),
that leave every element $\F(\O)$ globally invariant;
unbroken internal symmetries leave the vacuum state invariant.

Given a certain (compact) group $G$ of (unbroken) 
internal symmetries of $\F$,
the fixpoint net $\F^G$ defined by $\F^G(\O)=\F(\O)^G$ is 
an example of {\it subsystem}
(sometimes also called {\it subnet} or {\it subtheory} in the literature),
i.e. a net of (von Neumann) subalgebras of $\F$.
This is the typical situation allowing one to recover an
observable net from a field net via a principle of gauge invariance.
However, in certain situations one can easily produce examples of
subsystems
that can hardly be seen to arise in this way.
See e.g. the discussion in \cite{Wic,Ar,roberto1}.

In this work we address the problem of classifying subsystems
of a given net $\F$.
Some related work has been already done in
\cite{LaSc1,LaSc2,Dav,CDR,Car,seba2}.
Our main result states that 
if $\F$ satisfies certain structural properties 
then all the reasonably well-behaved subsystems 
morally
arise in the way explained above, 
namely they are fixpoints for a compact group action 
on $\F$ or on one component $\F_1$ 
in a tensor product decomposition 
$\F=\F_1 \otimes \F_2$.

We confine our discussion to nets $\F$ satisfying usual postulates
such as Poincar\'e covariance, Bisognano-Wichmann and the split property,
plus an additional condition,
the absence of nontrivial sectors,
whose meaning has been recently clarified in \cite{CDR}.
Our assumptions are sufficiently general to cover many
interesting situations,
including the well-known 
Bosonic free field models (massive or massless).
In particular
in the case of (finitely many) multiplets of the massive scalar free
fields we (re)obtain a classification result of Davidson \cite{Dav},
but with a different method of proof. 
Moreover our discussion applies to the massless case as well.   
In a different direction, we also 
provide a first solution to
a long-standing open problem, proposed by S. Doplicher, 
concerning the relationship between
an observable net $\A$ and
the subsystem $\C$
generated by the local energy-momentum tensor \cite{Dop,roberto1}.
As to the main ingredients, 
now $\A$ is required to have 
the split property and 
at most countably many superselection sectors,
all with finite statistical dimension 
\footnote{If one can rule out the occurrence of sectors with infinite
statistics for $\A$, the other two facts are easily implied by the split
property for the canonical field net $\F$, 
that is anyhow needed from the start to define the subsystem $\C$.}
(and Bosonic).

Still our assumptions 
are restrictive enough to rule out the occurrence of 
models with undesiderable features.
This allows us to overcome certain technical difficulties that
cannot be handled in too general (perhaps pathological) situations.

This paper is organized in the following way.
In the next section we 
describe our setup and 
collect some preliminaries.
The third section contains the stated classification result.
In the fourth section we present some applications. 
Some of the assumptions can be relaxed to some extent, 
at the price of much more complicated proofs and no sensible improvement.
We end the article with some brief comments and suggestions for future work.
An appendix is included to provide some technical results about 
scalar free field theories.

\section{Preliminaries}
Throughout this article we denote $\P$ the 
connected
component of the identity
of the Poincar\'e group in four spacetime dimensions and $\K$ the set of 
open double cones of $\RR^{4}$. We will denote the elements of $\P$ by pairs
$(\Lambda,x)$, where $\Lambda$ is an element of the restricted Lorentz 
group and $x\in \RR^4$ is a spacetime translation, or alternatively by a 
single letter $L$. Double cones and wedges will be denoted 
$\O$ and $\W$ respectively, with subscripts if necessary.
We consider a net $\F$ over $\K$, i.e. a map 
$\O \to \F(\O)$ from double cones to von Neumann algebras acting on a 
separable Hilbert space $\H$, satisfying the following  assumptions.
\begin{itemize}
\item[(i)] {\it Isotony.} If $\O_1 \subset \O_2,\;\O_1,\O_2 \in \K$, then 
\begin{equation}
\F(\O_1) \subset \F(\O_2). 
\end{equation}

\item[(ii)] {\it Locality.} If $\O_1,\O_2 \in \K$ and $\O_1$ is 
spacelike separated from $\O_2$ then 
\begin{equation}
\F(\O_1)\subset\F(\O_2)^\prime,
\end{equation}
 
\item[(iii)] {\it Covariance.} There is a strongly continuous unitary 
representation $U$ of $\P$ such that, for every $L\in \P$ and every 
$\O \in \K$, there holds 
\begin{equation}
U(L)\F(\O)U(L)^*=\F(L\O).
\end{equation}

\item[(iv)]{\it Existence and uniqueness of the vacuum.} There exists a 
unique (up to a phase) unit vector $\Omega$ which is invariant 
under the restriction of $U$ to the one-parameter subgroup of spacetime 
translations. 

\item[(v)]{\it Positivity of the energy.} The joint spectrum of the 
generators of the spacetime translations is contained in the closure 
$\overline{V}_+$ of the open forward light cone $V_+$. 

\item[(vi)] {\it Reeh-Schlieder property.}  The vacuum vector $\Omega$ is 
cyclic for $\F(\O)$ for every $\O \in \K$. 

\item[(vii)] {\it Haag duality.} For every double cone 
$\O \in\K$ there holds
\begin{equation}
\F(\O^\prime)=\F(\O)^\prime,
\end{equation}
where $\O^\prime$ is the interior of the spacelike complement of $\O$ and, 
for every open set ${\cal S} \subset \RR^{4}$, $\F({\cal S})$ denote the 
algebra defined by
\begin{equation}
\F({\cal S})=\vee_{\O \subset {\cal S}} \F(\O).
\end{equation}

\item[(viii)] {\it TCP covariance.} There exists an antiunitary involution 
$\Theta$ (the TCP operator) such that:
\begin{eqnarray}
\Theta U(\Lambda,x) \Theta=U(\Lambda,-x)\;\;\forall (\Lambda,x)\in \P;
\\
\Theta \F(\O) \Theta =\F(-\O).
\end{eqnarray}

\item[(ix)] {\it Bisognano-Wichmann property.} Let 
\[ \W_R=\{x\in \RR^4 : x^1>|x^0| \} \]
be the right wedge and let $\Delta$ and $J$ be the modular operator and
the 
modular conjugation of the algebra $\F(\W_R)$ with respect to $\Omega$, 
respectively. Then it holds:
\begin{eqnarray}
\Delta^{it}=U(\Lambda(- 2\pi t),0);
\\
J=\Theta U(R_1(\pi),0);
\end{eqnarray}
where $\Lambda(t)$ and $\R_1(\theta)$ are the one-parameter groups 
of Lorentz boosts in the $x^1$-direction and of spatial rotations around 
the first axis, respectively. 

\item[(x)] {\it Split property.} Let $\O_1,\O_2 \in\K$ be open double cones 
such that the closure of $\O_1$ is contained in $\O_2$ (as usual we write 
$\O_1\subset\subset \O_2$). Then there is a type I factor $\N(\O_1,\O_2)$ 
such that 
\begin{equation}
\F(\O_1)\subset\N(\O_1,\O_2)\subset \F(\O_2).
\end{equation}
\end{itemize} 

Using standard arguments (cf. \cite{DA}) it can be shown that the previous 
assumptions imply the irreducibility of the net $\F$, namely the algebra 
$\F(\RR^4)$ coincides with the algebra ${\rm B}(\H)$ of all bounded operators
on $\H$. Another easy consequence of the assumptions is that $\Omega$ is 
$U$-invariant. Moreover the algebra $\F(\W)$ is a factor 
(in fact a type ${\rm III}_1$ factor), for every wedge $\W$,
see e.g. \cite[Theorem 5.2.1]{Borch99}. 
Strictly speaking,
it is also possible to deduce (viii) 
from the other assumptions \cite[Theorem 2.10]{GuLo1}.

From Haag duality it follows that the algebra associated with a double 
cone coincides with intersection of the algebras associated to the wedges 
containing it, i.e.
\begin{equation}
\F(\O)=\cap_{\O\subset \W} \F(\W),
\end{equation}
for every $\O\in \K$. Thus our net $\F$ corresponds to a particular case 
of the AB-systems described in \cite{Wic}, see also \cite{TW}. Moreover 
the Bisognano-Wichmann property implies wedge duality, i.e.
\begin{equation} 
\F(\W)^{\prime}=\F(\W^{\prime}),
\end{equation}
for every wedge $\W$, where $\W^{\prime}$ denotes the interior of the 
causal complement of $\W$. 

Another important fact is that, due to the split property, the net $\F$
satisfies Property B for double cones: given 
$\O \subset\subset \tilde{\O}$, $\O,\tilde{\O} \in \K$, for each nonzero 
selfadjoint projection $E\in \F(\O)$ there exists an isometry 
$W \in \F(\tilde{\O})$ with $E=WW^*$. 
Moreover, for every nonempty open set 
${\cal S}\subset \RR^4$, the algebra $\F({\cal S})$ is properly infinite. 

\begin{definition} {\rm A {\it covariant subsystem} $\B$ of $\F$ is an 
isotonous net of nontrivial von Neumann algebras over $\K$, such that: 
\begin{eqnarray}
\B(\O)\subset\F(\O);
\\
U(L)\B(\O)U(L)^*=\B(L\O),
\end{eqnarray}
for every $\O\in\K$ and every $L\in \P$.}
\end{definition}
We use the notation $\B\subset \F$ to indicate that $\B$ is a covariant 
subsystem of $\F$.
As in the case of $\F$, for every open set ${\cal S} \subset \RR^4$ we
define
$\B({\cal S})$ by 
\begin{equation}
\B({\cal S})=\vee_{\O\subset {\cal S}} \B(\O).
\end{equation}
\begin{definition}{\rm We say that a covariant subsystem $\B$ of $\F$ 
is {\it Haag-dual} if
\begin{equation}
\B(\O)=\cap_{\O\subset \W} \B(\W)\;\;\forall \O\in \K.
\end{equation}
}
\end{definition}
If a covariant subsystem $\B$ is not Haag-dual, one can associate to it 
an Haag-dual covariant subsystem $\B^d$ (the {\it dual subsystem}) 
defined by
\begin{equation}
\B^d(\O)= \cap_{\O \subset \W} \B(\W),
\end{equation} 
cf. \cite{TW,Wic}. Note that $\B(\W)=\B^d(\W)$ for every wedge $\W$. 

Given a covariant subsystem $\B$ of $\F$ we denote $\H_\B$ the 
closure of $\B(\RR^4)\Omega$ and by $E_\B$ the corresponding orthogonal 
projection. It is trivial that the algebras $\B(\O),\;\O\in\K$ leave 
$\H_\B$ stable. Hence we can consider the reduced von Neumann 
algebras $\hat{\B}(\O):=\B(\O)_{E_\B},\;\O\in\K$ acting on the Hilbert space 
$\H_\B$ and denote $\hat{\B}$ the corresponding net.  
It is straightforward to verify that 
\begin{equation}
\B({\cal S})_{E_\B}=\vee_{\O\subset {\cal S}}\hat{\B}(\O),
\end{equation}
for every open set ${\cal S}\subset \RR^4$. Therefore the notation
$\hat{\B}({\cal S})$ 
is unambiguous. Moreover, due to the Reeh-Schlieder property (for 
$\F$), the map $B\in \B({\cal S}) \mapsto \hat{B}:=B_{E_\B} \in
\hat{\B}({\cal S})$, is an isomorphism of von Neumann algebras, 
whenever the interior ${\cal S}^\prime$ of 
the causal complement of $\cal S$ is nonempty. 

The following result is due in large part to Wichmann \cite{Wic} and 
 Thomas and Wichmann \cite{TW}.
 \begin{proposition} Let $\B$ be a Haag-dual subsystem of $\F$. Then 
 the following properties hold:
 \begin{itemize}
 \item[(a)] $\Theta$ and $U$ commute with $E_\B$. Accordingly we can 
 consider the reduced operators $\hat{\Theta}:=\Theta_{E_\B}$ and 
 $\hat{U}:=U_{E_\B}$ on $\H_\B$; 

\item[(b)] All the properties from (i) to (x) listed in the beginning of 
this section 
holds with $\F$, $\H$, $U$, $\Theta$, replaced by $\hat{\B}$, $\H_\B$, 
$\hat{U}$, $\hat{\Theta}$, respectively. 
\end{itemize} 
\label{eredita}
\end{proposition}
\noindent {\it Proof.} For (a) and (b), properties from (i) to (ix), we
refer the reader 
 to \cite{Wic} and \cite[Section 5]{TW}. Proving (x) for $\hat{\B}$ 
 corresponds to show that the split property is hereditary. This fact 
 is well known (cf. e.g. \cite[Section 5]{Do82}) but we include 
 here a proof for convenience of the reader. 

 Let $\O_1,\O_2\in \K$  be such that $\O_1\subset\subset\O_2$. It is 
 sufficient to show that there is a faithful normal product state on 
 $\hat{\B}(\O_1)\vee\hat{\B}(\O_2)^\prime$, i.e. a faithful normal state 
 $\phi$ satisfying
 \begin{equation}
 \phi(B B^\prime)=\phi(B)\phi(B^\prime) \;\;\forall B\in \hat{\B}(\O_1),\;
 \forall B^\prime \in \hat{\B}(\O_2)^\prime,
 \end{equation}
 see e.g. \cite{DoLo84}. 
 $\hat{\B}$ satisfies Haag duality and 
 \[\hat{\B}(\O_1)\vee\hat{\B}(\O_2^\prime)=
 [\B(\O_1)\vee\B(\O_2^\prime)]_{E_\B}\]
 is isomorphic to $\B(\O_1)\vee\B(\O_2^\prime)$, 
 being $\H_\B$ separating for the latter algebra.
 Therefore it remains to show the existence of a 
 faithful normal product state on $\B(\O_1)\vee\B(\O_2^\prime)$. 
 This trivially follows 
 from the existence of a faithful normal product state for 
 $\F(\O_1)\vee\F(\O_2^\prime)$, which is a consequence of the split 
 property for $\F$. 
\footnote{A similar argument shows that split for wedges 
(cf. \cite{Mu2}) 
is inherited by subsystems satisfying wedge duality;
here the space-time dimension is not important.}
\qed

From the previous proposition it follows that if $\B$ is Haag-dual then 
$\hat{\B}$ satisfies Haag duality.\footnote{This is not true in two
spacetime dimensions.} 
It is quite easy to show that also the converse is true. This 
remark should make it clear that considering only Haag-dual subsystems is 
not a too serious restriction. 

If $\B$ is a covariant subsystem of $\F$, we can consider the net $\B^c$
defined by
\begin{equation}
\B^c(\O)=\B(\RR^4)^\prime \cap \F(\O),
\end{equation}
cf. \cite{Dav,Borch99}.
If $\B^c$ is trivial, then we say that $\B$ is {\it full} (in $\F$).
If $\B^c$ is nontrivial, then it is easy to check that it is a Haag-dual
covariant subsystem of $\F$ ({\it the coset subsystem}).
It follows from the definition that $\B \subset \B^{cc}$, 
and $\B^c = \B^{ccc}$.

For later use it is convenient to introduce the notions 
of tensor product and of unitary equivalence of two
nets.
Let $\F_1$ and $\F_2$ be two nets 
acting on $\H_1$ and $\H_2$ respectively, 
and let $U_1, U_2$ and $\Omega_1, \Omega_2$ the corresponding
representations of the Poincar\'e group and the vacuum vectors.
By {\it tensor product} of nets $\F_1 \otimes \F_2$ we mean the net
$\K\ni\O\mapsto\F_1(\O)\otimes\F_2(\O)$ acting on $\H_1 \otimes \H_2$
together with the representation $U_1 \otimes U_2$ of $\P$
and the vacuum $\Omega_1\otimes\Omega_2$. 
It follows that $\F_1\otimes\F_2$ 
satisfies properties (i)--(x) if $\F_1$ and $\F_2$ do so.
We say that $\F_1$ and $\F_2$ are {\it unitarily equivalent} if 
there exists a unitary operator $W: \H_1 \to \H_2$ with
$W\F_1(\O)W^*=\F_2(\O) \ (\O \in \K)$, $WU_1(L)W^*=U_2(L)$.
Note that since the vacuum is unique up to a phase, one can always choose
$W$ so that $W\Omega_1=\Omega_2$.

\section{General Classification Results}
\label{gen}
In this section we consider a net $\F$ satisfying all the properties 
(i)--(x) described in the previous section. Moreover we will assume 
the following condition (cf. \cite{CDR}):
\begin{itemize}
\item[(A)] Every representation of (the quasi-local $C^*$-algebra) $\F$
satisfying the DHR selection criterion 
is a multiple of the vacuum representation.\footnote{For the basic
notions concerning 
the DHR theory of superselection sectors we refer the reader to \cite{Ha}
and references therein.} 
\end{itemize}

Let us observe that condition (A) is equivalent 
to the seemingly weaker condition that all the irreducible
representations satisfying the selection criterion are 
equivalent to the vacuum representation.
This is a consequence of the fact that 
the irreducible representations occurring
in the direct integral 
decomposition of a localized\footnote{In this article the word 
{\it localized} referred to representations or endomorphisms means  
{\it localized in double cones}.}
representation are localized a.e.
(see \cite[Appendix B]{KLM}).

Now let $\B$ be a Haag-dual, covariant subsystem of $\F$ and let $\pi$ 
be the corresponding representation of $\hat{\B}$ in $\H$, i.e. 
the representation defined by $\pi(\hat{B}) =B$ for 
$B\in \cup_{\O \in \K} \B(\O)$. We denote 
$\pi^0$ the identical (vacuum) representation of $\F$ on $\H$ and 
$\pi_0$ the vacuum representation of $\hat{\B}$,
i.e. its identical representation on $\H_\B$. 
The following result is already known (see e.g. \cite{CDR}) but we include
a proof for the sake of completeness.

\begin{lemma} $\pi$ satisfies the DHR criterion.
\label{DHRlemma}
\end{lemma}
\noindent{\it Proof.} For every $\O \in \K$ the von Neumann algebras 
$\B(\O^\prime)$ and $\hat{\B}(\O^\prime)$ are isomorphic. Moreover, as 
noted in the previous section, these von Neumann algebras are properly 
infinite with properly infinite commutants. By \cite[Theorem 7.2.9.]{KR2} 
and \cite[Proposition 9.1.6.]{KR2} we can find a unitary operator 
\[U_\O: \H_\B \to \H \]
such that 
\[ U_\O \hat{B} {U_\O}^* =B\;\; \forall B \in \B(\O^\prime). \] 
Hence if $\O_1\in \K$ is contained in $\O^\prime$ there holds 
\[ \pi_0 (\hat{B})= {U_\O}^* \pi(\hat{B}) {U_\O} \;\;\forall \hat{B}\in 
\hat{\B}(\O_1).\]
Actually, this is the DHR criterion.\qed

\begin{proposition} For every irreducible localized transportable morphism 
$\sigma$ of $\hat{\B}$, $\pi_0 \circ \sigma$ is equivalent to a 
subrepresentation of $\pi$. 
Moreover $\sigma$ is covariant with positive energy and 
it has finite statistical dimension. 
\label{hearth}
\end{proposition}

\noindent{\it Proof.} Since $\pi$ satisfies the DHR criterion we can 
find a transportable localized morphism $\rho$ of $\hat{\B}$ such that
there holds the unitary equivalence 
\begin{equation}
\pi\simeq \pi_0 \circ \rho,
\label{pi}
\end{equation}
cf. \cite[Proposition 3.4.]{LR}. 

Let us consider the extension $\hat{\sigma}$ of $\sigma$ to $\F$ 
\cite{CDR}, cf. \cite{LR}. 
Then the assumption (A) for $\F$ imply that 

\begin{equation}
\pi^0\circ\hat{\sigma}\simeq\oplus_i \pi^0,
\label{sigmaext}
\end{equation}
where the index $i$ in the direct sum on the r.h.s.
runs over a set whose
cardinality is at most countable. 
Restricting these representations to $\B$ we find 
\begin{equation}
\pi\circ\sigma\simeq\oplus_i\pi
\end{equation}
and therefore using equation \ref{pi} 
\begin{equation}
\rho \sigma \simeq \oplus_i\rho.
\label{rhosigma}
\end{equation}   

Since $\rho$ contains the identity sector we have 
$\sigma \prec \rho \sigma$ and hence 
\begin{equation}
\sigma \prec \oplus_i \rho.
\label{contained}
\end{equation} 

Thus, being $\sigma$
arbitrary, every irreducible representation of $\hat{\B}$ satisfying the 
DHR criterion is contained in a countable multiple of $\rho$. The latter 
multiple is a representation on a separable Hilbert space. Hence there are
at most countably many irreducible sectors of $\hat\B$.

Being $\pi$ a direct integral of irreducible DHR representations
\cite[Appendix B]{KLM} and appealing to some standard 
arguments (see e.g. \cite{Dix1,Dix2}) one gets that $\pi$ is in fact a
direct sum. 
From 
equation \ref{contained} it is not difficult to show that,
being $\sigma$ irreducible, we have $\sigma\prec \rho$ 
i.e. $\pi_0 \circ\sigma$ is unitarily equivalent to a subrepresentation
of $\pi$.

Since $\B$ is covariant $\pi$ is covariant with positive energy. 
We have to show that every 
irreducible
subrepresentation has the same property,
cf. \cite{BCL}.
 Since the action induced by
the representation $U$ of the Poincar\'e group leaves 
$\B(\RR^4)$ globally invariant it leaves globally invariant also its 
centre. 
Being the latter purely atomic
(due to the decomposition of $\pi$ into irreducibles) 
and $\P$ connected,
it follows that the orthogonal projection 
$E_{[\sigma]}\in \B(\RR^4)' \cap \B(\RR^4)$ onto the 
isotypic subspace corresponding to $\sigma$ must commute with $U$. 
Let $U_{[\sigma]}$ and $\pi_{[\sigma]}$ be the restrictions to 
$E_{[\sigma]} \H$ of $U$ and $\pi$ respectively. Then we have the unitary 
equivalence
\begin{equation}
\pi_{[\sigma]}\simeq (\pi_0 \circ \sigma) \otimes I \ .
\end{equation}
Moreover, using the relation
\begin{equation}
U_{[\sigma]}(L)\pi_{[\sigma]}(\hat{B})U_{[\sigma]}(L)^*
=\pi_{[\sigma]}(\hat{U}(L)\hat{B}\hat{U}(L)^*),
\end{equation}
where $B \in \cup_{\O\in\K} \B(\O), \; L \in \P$,
and a classical result by Wigner on projective unitary
representations of $\P$ \cite{Wig,Barg},
it is quite easy to show that
\begin{equation}
U_{[\sigma]}(L)\simeq
U_\sigma(L)\otimes X_\sigma (L),
\end{equation}
where $U_\sigma$ and $X_\sigma$ are unitary continuous representations 
of (the covering group of) $\P$ 
and $U_\sigma$ is such that 
\begin{equation}
U_\sigma(L)\sigma(\hat{B})U_\sigma(L)^*
=\sigma(\hat{U}(L)\hat{B}\hat{U}(L)^*).
\end{equation}
Since $U_{[\sigma]}$ satisfies the spectral condition, both 
$ U_\sigma$ and $X_\sigma$ have to satisfy it.\footnote{This follows from
the fact that if $S_1$ and $S_2$ are two orbits of the restricted Lorentz 
group such that $S_1+S_2\subset \overline{V}_+$ then
$S_1\subset\overline{V}_+$
and $S_2\subset \overline{V}_+$.}
Hence $\sigma$ is covariant with positive energy.  

Finally, from $\rho \sigma\simeq \sigma \rho$ and equation \ref{rhosigma} 
it follows that ${\rm id}\prec \sigma \rho$.
Therefore, being $\sigma$ covariant with positive energy, 
it has finite statistical dimension because of \cite[prop. A.2]{DHR4}. 
\qed

A related result has been independently obtained by R. Longo, 
in the context of nets of subfactors \cite{longo}.

\medskip
Let $\F_\B$ be the canonical field net of $\hat{\B}$ as defined in
\cite[Section 3]{DoRo90}. 
In natural way $\F_\B$ can be considered as a Haag-dual 
subsystem of $\F$ containing $\B$ \cite[Theorem 3.5]{CDR}. 
In fact one finds that 
$\F_\B(\O)$ coincides with the von Neumann algebra generated by the 
family of Hilbert spaces $\H_{\hat{\sigma}}$ in $\F$, where $\sigma$ runs
over all the transportable morphisms of $\B$ which are localized in $\O$
and $\hat{\sigma}$ denotes the functorial extension of $\sigma$ to $\F$.
From the fact that the latter extension commutes with spacetime symmetries, 
namely $(\sigma_L)\hat{} = (\hat{\sigma})_L $ 
for every $L \in \P$ it is also easy to show that $\F_\B$ is a covariant 
subsystem.
(Besides, by \cite[Proposition 2.1]{roberto}
$\F_\B$ coincides with its covariant
companion, cf. \cite{DoRo90}.)

\begin{theorem} $\hat{\F_\B}$ has no irreducible DHR sectors other than
the vacuum. 
\end{theorem}

\noindent{\it Proof.} 
By the previous proposition it is enough to consider sectors
with finite statistical dimension.
Let $\R$ be the canonical field algebra of  $\hat{\F_\B}$. 
Then $\R$ is a Haag-dual covariant subsystem of $\F$,  
and as such it inherits the split property.
By the results discussed in \cite{ext} 
this is sufficient\footnote{This idea is not new,
see e.g. \cite[Section 2]{Mu}, 
however some technical difficulties are circumvented
when the assumptions made in this paper are used.}
to deduce that $\F_\B=\R$.\footnote{Alternatively, the same
result may be deduced combining Proposition \ref{hearth} with 
\cite{CDR}.}
In fact the group $\tilde{G}$ of 
the (unbroken) symmetries of $\R$
extending the
gauge automorphisms of $\F_\B$ is compact 
in the strong operator topology
by (the proof of) \cite[Theorem 10.4]{DoLo84},
and obviously $\R^{\tilde{G}}=\B$. 
The conclusion follows
by the uniqueness of the canonical field net \cite{DoRo90}.
\qed

\begin{theorem} 
There exists a unitary isomorphism of $\F$ with 
$\hat{\F_\B} \otimes \hat{\B^c}$.
In particular $\F_\B=\B^{cc}$, and 
if $\B$ is full\footnote{Irreducible subsystems, 
namely those satisfying 
$\B' \cap \F={\mathbb C}$, are full.} 
in $\F$ then $\F_\B=\F$.
\label{splitting} 
\end{theorem}

\noindent{\it Proof.} Let $\tilde\pi$ be the 
representation
of $\hat{\F_\B}$ on $\H$ (the vacuum Hilbert space of $\F$)
arising from the embedding $\F_\B \subset \F$ 
and $\tilde{\pi_0}$ the vacuum representation of $\hat{\F_\B}$ on
$\H_{\F_\B}\subset \H$. 
By the previous theorem $\hat{\F_\B}$ has no nontrivial sectors. Moreover  
Lemma \ref{DHRlemma} applied to $\F_\B$ instead of $\B$ implies that  
$\tilde{\pi}$ is (spatially) equivalent
to a multiple of $\tilde{\pi_0}$ 
and therefore to $\tilde{\pi_0} \otimes {\rm I}$, on 
$\H_{\F_\B} \otimes \H_1$, where $\H_1$ is a suitable Hilbert space. 
Let $W: \H \to \H_{\F_\B} \otimes \H_1$ be a unitary
operator implementing this equivalence.
For every double cone $\O$ there holds 
\begin{equation}
\hat{\F_\B} (\O^{\prime}) \otimes {\rm I} \subset  
\tilde\F(\O^{\prime} )
\end{equation}
where $\tilde{\F}(\O)=W\F(\O)W^*$.
Therefore, using Haag duality for $\hat{\F_\B}$, 
\begin{equation}
\hat{\F_\B} (\O) \otimes {\rm I} \subset  \tilde\F (\O)\subset 
\hat{\F_\B} (\O) \otimes {\rm B}(\H_1 ). 
\end{equation}
It follows that 
\begin{equation}
\hat{\F_\B} (\W) \otimes {\rm I} \subset  \tilde\F (\W)\subset 
\hat{\F_\B} (\W) \otimes {\rm B}(\H_1 ). 
\end{equation}

The algebras of wedges are factors. By the results in \cite{GeKa} 
(cf. also \cite{StZs})
there exists a von Neumann algebra  
$\M (\W)\subset {\rm B}(\H_1 )$ such that
\begin{equation}
\tilde\F (\W) =\hat{\F_\B} (\W) \otimes \M (\W).
\end{equation}
Taking on both sides of this equality the intersection over all the wedges
containing a given $\O \in \K$ we find
\begin{equation}
\tilde\F (\O) =\hat{\F_\B} (\O) \otimes \M (\O),
\end{equation}
where 
\begin{equation}
\M(\O)=\cap_{\O\subset \W}\M(\W). 
\end{equation}
Now, using the commutant theorem for von Neumann tensor products, 
it is straightforward to show that 
\[{\rm I}\otimes \M(\O)=W\B^c (\O)W^*\]
for every $\O \in \K$.
The previous equation implies the existence of a representation $\tau$ of 
$\hat{\B^c}$ on $\H_1$ such that $WBW^*={\rm I}\otimes \tau (\hat{B}), \
B\in \B^c(\O)$ for every $\O\in\K$. Moreover, since $\M$ acts irreducibly
on $\H_1$ and the vacuum representation $\pi^c$ of $\hat{\B^c}$ is
contained in ${\rm I}\otimes \tau$, $\tau$ is spatially isomorphic to
$\pi^c$ and thus the mapping  $\O \to \M(\O)$ gives a net unitarily
equivalent to $\hat{\B^c}$. Therefore without loss of generality we can
assume that $\H_1=\H_{\B^c}$ and that $W\F(\O)W^*=\hat{\F_\B}(\O)\otimes
\hat{\B^c}(\O), \; \O\in \K$.
The conclusion follows noticing that $WUW^*=U_{E_{\F_\B}} \otimes
U_{E_{\B^c}}$. Here we omit the easy details.
\qed

Applying the previous theorem to $\B^c$ in place of $\B$ we get that
$\B^c$ as no nontrivial sectors, since $\F_{\B^{c}}=\B^{ccc}=\B^c$.

\begin{corollary} Let $\B$ be a Haag-dual covariant subsystem of $\F$,
then the net of inclusions $\K \ni \O \mapsto \B(\O) \subset \F(\O)$
is (spatially) isomorphic to 
$\O \mapsto \hat{\F_\B}(\O)^G \otimes I \subset \hat{\F_\B}(\O)\otimes
\hat{\B^c}(\O)$, where $G$ is the canonical gauge group of $\hat\B$.
\end{corollary}

\begin{corollary} If $\B$ is a Haag-dual covariant subsystem of $\F$
and if $\F_\B $ is full (in particular if $\B$ is full) 
then there exists a compact group $G$ 
of unbroken internal symmetries of $\F$ such that
$\B=\F^{G}.$
\label{classification}
\end{corollary}

Now let $\C$ be the 
(local) 
net generated by the canonical implementations of the translations
on $\F$ \cite{roberto1}. It is a covariant subsystem of $\F$. 
Since $\C$ is 
(irreducible thus)
full in $\F$ and $\C^d\subset\F^{G_{\rm max}}$, where
$G_{\rm max}$ is the (compact) group of all unbroken internal symmetries
of $\F$,  we have

\begin{corollary} In the situation described above it holds
\begin{equation}
\C^d=\F^{G_{\rm max}}.
\end{equation}
\end{corollary}

\section{Applications}

\subsection{Free fields}
Our standing assumptions are satisfied in the case where $\F$ 
is generated by a finite set
of free scalar fields \cite{Dr,BDLR} 
and also by suitable infinite sets of such fields \cite{DP}. 
They are also satisfied in other Bosonic theories, 
e.g. when $\F$ is generated by the free electromagnetic field, 
see \cite{BDLR}.

Therefore from our Corollary \ref{classification} 
one can obtain all the results in \cite{Dav} in the case of
full subsystems, even without assuming the existence of a mass gap.
Concerning subsystems that are not full, 
one has to study the possible decompositions 
\begin{equation}
\hat{\F_\B}(\O)\otimes \hat{\B^c}(\O) =\F(\O)
\end{equation}
(up to unitary equivalence).
In the case where $\F$ is generated by a finite set of free scalar fields,
it turns out that $\F_\B$ and $\B^c$ are always free scalar theories 
generated by two suitable disjoint subsets of the generating
fields of $\F$.
We present a detailed proof of this fact in the appendix.
\footnote{Davidson obtained this result
in the purely massive case \cite{Dav}.}
In particular, if $\F$ is generated by 
a single scalar free field $\varphi(x)$ of mass $m \geq 0$, 
no such nontrivial decomposition is possible and hence
all of the
subsystems of $\F$ are full. 
Accordingly, in this case, the unique Haag-dual covariant proper subsystem
of $\F$ is the fixed point net $\F^{\ZZ_2}$ under the action of the
group of (unbroken) internal symmetries.

\medskip
Note that when $m=0$
there are covariant subsystems which are not Haag-dual.
For instance the subsystem $\A \subset \F$ generated 
by the derivatives $\partial_\mu\varphi(x)$ is Poincar\'e 
covariant but not Haag-dual and in fact one has 
$\F=\A^d$ \cite{BDLR}.
However it can been shown that conformally covariant subsystems of $\F$
are always Haag-dual. Actually the latter fact still holds in a more
general context.

\subsection{Theories with countably many sectors}
Summing up,
we have shown a classification result for 
Haag-dual
subnets of a 
purely Bosonic net with trivial superselection structure 
(including infinite statistics)
and with the split property.
Moreover we have exhibited an important class of examples, 
namely (multiplets of) the free fields, 
to which our results apply.
This is already quite satisfactory.
One can consider a more general situation in which $\F$ is the
canonical field net of an observable net $\A$.
A closely related problem is, of course, to look for the structural
hypotheses on $\A$
ensuring that $\F=\F_\A$ will have the required properties.
It has been known for some time 
that if $\A$
has only a finite number of irreducible DHR sectors
with finite statistical dimension (i.e. $\A$ is rational), 
all of which are Bosonic,
then $\F$ (is local and) has no nontrivial DHR sectors with finite
statistical dimension \cite{roberto,Mu}. 
This result is not sufficient for our purposes, because it does not rule
out the possible presence of irreducible DHR representations of $\F$
with infinite statistical dimension.
However, a solution to this problem can be achieved by using the stronger
results given in \cite{CDR}.

\begin{theorem}
Let $\A$ be a local net satisfying the split property and Haag duality in
its (irreducible) vacuum representation.
If $\A$ has at most countably many irreducible (DHR) sectors, 
all of which are 
Bosonic and 
with finite statistical dimension, 
then any 
sector of $\A$ is a direct sum of irreducible 
sectors.
Moreover, the canonical field net $\F$ of $\A$ has no nontrivial 
sectors with {\it any} (finite or infinite) statistical dimension.
\end{theorem}

\noindent {\it Proof.} 
In view of \cite[Theorem 4.7]{CDR} it is enough to show the first statement. 
But using the split property and the bound on the number of inequivalent
sectors, this follows arguing as in the proof of Proposition \ref{hearth}.
\qed

\medskip
This result 
\footnote{As in \cite{roberto}, in the case of rational theories
a different argument could be given 
when the local algebras are factors, 
based on a restriction-extension argument 
(cf. \cite[Lemma 27]{KLM}).}
shows that $\F$ satisfies the condition (A) of section
\ref{gen}.
Moreover if $\A$ satisfies all of the conditions (i)-(vii)
then the same is true for $\F$ \cite{DoRo90}.
In order to apply the above result about classification
of subsystems and solve the problem about local charges, 
we need to know conditions on $\A$ implying the validity of properties
(viii)-(x).
Concerning (x), 
it would be a consequence of the split property for $\A$ 
if $G$ were finite and abelian.
In other cases one can invoke some version of nuclearity for $\A$,
implying that $\F$ is split \cite{BuDA}.
But it is also necessary to know if the existence of a TCP symmetry 
and the special condition of duality for $\A$ imply the same for
its canonical field system $\F$.
The relationship between the validity of conditions (viii)-(ix) for $\A$
and its canonical field system $\F$ has been 
discussed in \cite{Kuc,Kuc2} 
(the TCP symmetry has been also treated
in \cite{CoDA} under milder hypotheses).
The conclusion is that if $\A$ satisfies the usual axioms
(and all its sectors are covariant), 
moreover it is purely Bosonic 
and satisfies a suitable version of nuclearity
(implying, among other things, the existence of at most countably many
sectors), 
TCP covariance and the Bisognano-Wichmann property, 
then we know how to classify all the subsystems of $\F$ satisfying Haag
duality.

\begin{corollary} Let $\A$ be an observable net 
satisfying the properties (i)-(ix) above,
without DHR sectors with infinite statistical dimension 
or para-Fermi statistics of any finite order,
whose 
(Bosonic)
canonical field net $\F$ has the split property.
Then, if $\C$ is the net generated by the local energy-momentum tensor,
one has
$$\C^d=\F^{G_{\rm max}}.$$
Moreover
$\A=\C^d$ 
if and only if 
$\A$ has no proper full 
Haag-dual
subsystem
(in which case $\A$ has no unbroken internal symmetries).
\end{corollary}

\noindent{\it Proof.}
Since $\A$ satisfies the split property and has at most countably 
many sectors, all with finite statistics,
the first statement follows by the previous result and Corollary 3.3.
If $G$ denotes the canonical gauge group of $\A$,
so that $\A=\F^G$,
the equality $\A=\C^d$ is equivalent to the equality $G=G_{\rm max}$,
which, due to Corollary \ref{classification}, means that 
there is no proper subsystem of $\A$ full (or irreducible) in $\F$.
To complete the proof we only need to show that every full subsystem of 
$\A$ is full in $\F$, when $G=G_{\rm max}$. Let $\B$ be a (Haag-dual)
subsystem of $\A$. Due to the results in the previous section, for every 
wedge $\W$ the inclusions
\[ \B(\W)\subset\A(\W)\subset \F(\W) \]
are spatially isomorphic to
\[ \hat{\B}(\W)\otimes I \subset \tilde{\A}(\W) \subset \hat{\F_\B}(\W)
\otimes \hat{\B^c}(\W), \]
with $\tilde\A$ isomorphic to $\A$.
Moreover, from $G=G_{\rm max}$ it follows that
\[ \tilde{\A}(\W) \subset \hat{\B}(\W)\otimes \hat{\B^c}(\W). \]
Arguing as in the proof of theorem \ref{splitting} we find that if 
$\B$ is not full in $\F$ then for every $\O\in \K$, the algebra
$\B(\RR^4)' \cap \A(\O)$ is nontrivial. It follows that $\B$ is not 
full in $\A$.  
\qed

\section{Comments on the assumptions}
Some of the results of the previous sections are in fact still true
even after relaxing some conditions.
We will briefly discuss some aspects here.

\medskip
The hypothesis (x) is useful to derive property B 
(also for the subsystems),
to apply the results in \cite{KLM}
and also to define the local charges. 
If we renounce to 
(x), and possibly (A),
taking $\F$ as the DHR field algebra of $\A \supset \B$ in its vacuum
representation on $\H$
(here it is not even essential to require the condition of covariance, 
nor the additional assumptions of the main theorem in \cite{CDR}),
it is still possible to deduce that $\tilde\pi \simeq \tilde{\pi}_0
\otimes I$ 
as in the proof of Theorem \ref{splitting}. 
For this purpose one needs to know that $\A$ and $\B$ both satisfy
property B,
and that $\tilde\pi$ 
in restriction to $\B$ 
(thought of as a representation of $\hat\B$)
is quasi-contained in the canonical
embedding of $\hat\B$ into its field net.
By the results in \cite{CDR},
the latter property holds 
if it is possible to rule out the occurrence of
representations with
infinite statistics for 
$\hat\B$ acting on $\H$
(e.g. if $[\A:\B]<\infty$ in the case of nets of subfactors).
In fact we don't even need to know a priori 
that $\pi$ satisfies the DHR selection criterion.
Relaxing covariance is necessary to discuss QFT on 
(globally hyperbolic) curved spacetimes.
Possibly
results resembling those
presented here should 
hold also in that context (cf. \cite{GLRV}).

\medskip
The Bisognano-Wichmann property for $\F$ 
and TCP covariance
may also be relaxed,
but,
for the time being,
$\F$ and the considered subsystems always have to satisfy Haag duality in
order to deduce some nice classification result.

\medskip
However,
let us discuss the inheritance of the split property in a slightly more
general situation.
We start with a 
subsystem $\B \subset \F$, 
but now both $\F$ and $\B$ are only assumed to satisfy essential
duality (cf. \cite{Ha})
in their respective vacuum representation, 
namely $\F^{d}=\F^{dd}$ and
$(\hat{\B})^{d}=(\hat{\B})^{dd}$ 
(this is consistent with the notation adopted in the previous sections).
Moreover we require the split property for $\F^d$.
In the situation where one has an embedding 
of $(\hat{\B})^{d}$ inside $\F^d$, 
\footnote{This may be true or not and is related to the
validity of the equality $(\hat{\B})^{d}=(\B^d)\hat{}$.}
we may  deduce the split property for $(\hat{\B})^{d}$
by our previous argument.
For instance if $\F$ satisfies the Bisognano-Wichmann property 
(thus in particular wedge duality, which implies essential duality),
then $\hat\B$ satisfies the same property as well \cite{Wic}
and moreover there exists
the embedding alluded above,
therefore the split property for $\F^d$
entails the split property for $(\hat{\B})^{d}$.\footnote{As a matter 
of fact, 
the same argument goes through when we just have essential duality
for $\F$ and wedge duality for $\hat\B$, 
see e.g. 
\cite[Section 3]{CDR}.}

\section{Outlooks}
In this article we have not discussed graded local (Fermionic) nets.
As far as we can see, it should be possible to obtain classification
results also in this case, 
once the natural changes 
in the assumptions, the statements and the proofs
are carried out.

In the situation described in the present paper the index of a subsystem
is clearly always infinite, or an integer.
Moreover any integer value is in fact realized\footnote{To see this,
consider the fixpoint net of the complex scalar free field
under the subgroup ${\mathbb Z}_n$ 
of the gauge group $S^1$.}.
In a broader context 
(e.g. inclusions of conformal nets on $S^1$),
the computation of the set of possible index values for
subsystems
seems an interesting
problem. 
In the case of concrete models many calculations are now available.
We hope to return on these subjects in the future.

\appendix
\section{Appendix}

In this appendix we study the possible tensor product decompositions of a 
net generated by a finite number of scalar free fields.

We consider a net $\O \mapsto \F(\O)$, 
acting irreducibly on its vacuum Hilbert space $\H$,
generated by a finite family of Hermitian 
scalar free fields 
$\varphi_1(x)$, $\varphi_2(x)\ldots$, $\varphi_n(x)$,
where 
$n=n_1+n_2+ \ldots +n_k$ and
$\varphi_1(x), \ldots$, $\varphi_{n_1}(x)$ have mass $m_1$,
$\varphi_{n_1+1}(x), \ldots$, $\varphi_{n_1+n_2}(x)$ have mass $m_2$,
and so forth, 
and $0 \leq m_1 < \ldots < m_k$.

Accordingly, for each $\O\in\K$, $\F(\O)$ is the von Neumann algebra
generated by the Weyl unitaries 
$e^{i \varphi_j(f)}$ for $j=1,\ldots,n$
and real-valued $f\in {\cal S}(\RR^4)$ with support in $\O$.

We denote $U, \Theta, \Omega$ the corresponding representation of $\P$,
TCP operator and vacuum vector respectively.

For every $i$ we let ${\k}_i$ be the closed subspace of $\H$ 
generated by the vectors $\varphi_i(f)\Omega$ 
with $f \in {\cal S}(\RR^4)$.

Each $\k_i$ is $U$-invariant, 
and the restriction $V_i$ of $U$ to $\H_i$
is the irreducible representation of $\P$ with spin $0$ and corresponding
mass. 

Moreover the generating fields 
are chosen so that $\k_i$ is orthogonal to $\k_j$ for $i \neq j$.

If $\k=\oplus_{i=1}^n \k_i$ and $V=\oplus_{i=1}^n V_i$, then 
$\H$ can be identified with the (symmetric) Fock space $\Gamma(\k)$
and $U$ with the second quantization representation $\Gamma(V)$,
see e.g. \cite{RS2}.

If $\F_i$ is the covariant subsystem of $\F$ generated by $\varphi_i(x)$,
then $\H_{\F_i}$ can be identified with $\Gamma(\k_i)$ and 
from the relation $\F(\O)=\vee_i \F_i(\O)$ and the properties of the
second quantization functor it follows that the net $\F$ is isomorphic
to $\hat{\F_1} \otimes \ldots \otimes \hat{\F_n}$ on 
$\otimes_i \Gamma(\k_i)$.

Note that there is some freedom in the choice of the generating fields,
corresponding to the internal symmetry group 
$G={\rm O}(n_1) \times \ldots \times {\rm O}(n_k)$.

Let $E_{m_h}$ be the orthogonal projection from $\H$ onto 
$\k_{m_h}:=\oplus_{i=n_{h-1}+1}^{n_{h-1}+n_h}\k_i$,
where by convention $n_0=0$.
For each $m \geq 0$, let $P_m$ be the orthogonal projection onto 
${\rm Ker}(P^2 - m^2)$, where $P^2$ 
denotes the mass operator corresponding to $U$.
It is not difficult to see that $P_m(\k + \CC\Omega)^\perp=0$ 
by a direct 
calculation on the $k$-particles subspaces of $\H$
(note that $P_m=0$ whenever $m\notin \{0\}\cup \{m_1,\ldots,m_k\}$).
It follows that
$P_{m_h}=E_{m_h}$ if $m_h>0$, 
while for $m_h=0$ we have $P_{m_h}=E_{m_h} + P_\Omega$
where $P_\Omega \in U(\P)'\cap U(\P)''$ 
is the orthogonal projection onto $\CC\Omega$.
In particular, 
for any $h\in\{1,\ldots,k\}$ we have
$E_{m_h} \in U(\P)'\cap U(\P)'' .$

The following simple lemma will be used to study the tensor product
decomposition of $\F$.

\begin{lemma}
Let $U_1$ and $U_2$ be subrepresentations of $U$ on subspaces $\H_1$ and
$\H_2$ of $\H$ both orthogonal to $\CC\Omega$. Then there are no
eigenvectors for the mass operator corresponding to the representation
$U_1 \otimes U_2$.
\label{nomass}
\end{lemma}

\noindent{\it Proof.}
We consider the net $\tilde\F=\F \otimes \F$ and the corresponding
representation $\tilde{U}=U\otimes U$ of $\P$.
Obviously the net $\tilde\F$ is of the same type as $\F$, with the same
masses but different multiplicities.
$U_1 \otimes U_2$ is a subrepresentation of $\tilde{U}$ on $\H_1
\otimes\H_2$.
If $\tilde{P}^2$ is the mass operator corresponding to $\tilde{U}$
and $\tilde{P}_m$ is the orthogonal projection onto 
${\rm Ker}(\tilde{P}^2- m^2)$, we only have to show that
for every $m\geq 0$ we have $\tilde{P}_m \H_1 \otimes \H_2 =0$.
But this follows by the discussion in the last paragraph before the
statement, since $\H_1 \otimes \H_2$ is orthogonal to
$\CC(\Omega\otimes\Omega)+\tilde{\k}$
where $\tilde{\k}=\k\otimes \Omega +\Omega\otimes \k$
is the one-particle subspace of $\H\otimes\H$.
\qed

\medskip
We are now ready to study the possible tensor product decompositions
$\F_A \otimes \F_B$ of $\F$.
In the sequel we assume to have such a decomposition, and 
deduce some consequences.

Then $\H$ is given by $\H_A \otimes\H_B$
so that $\Omega=\Omega_A \otimes \Omega_B$
and $U=U_A \otimes U_B$.

We set $\H_A={\mathbb C}\Omega_A \oplus \tilde\H_A$, 
and analogously for $\B$, 
so that 
$\H=\H_A\otimes\H_B =
\CC \Omega \oplus (\Omega_A \otimes \tilde\H_B) \oplus 
(\tilde\H_A \otimes\Omega_B) \oplus (\tilde\H_A \otimes \tilde\H_B)$.
We also set $F_0=P_\Omega$, $F_A=[\Omega_A\otimes\tilde\H_B]$,
$F_B=[\tilde\H_A\otimes\Omega_B]$, $F_{AB}=[\tilde\H_A \otimes
\tilde\H_B]$.
Notice that these orthogonal projections commute not only with $U$
but also with $\Theta$.

\begin{lemma} 
For each $h=1,\ldots,k$  it holds  $E_{m_h}F_{AB}=0$.
\end{lemma}
\noindent{\it Proof.} It is an immediate consequence of Lemma
\ref{nomass}.\qed

\medskip

Since $E_{m_h}F_0=0$, the previous lemma implies that
$E_{m_h}(F_A+F_B)=E_{m_h}$, for $h=1,\ldots ,k$. 
This amounts to say that  
$\k\subset \tilde{\H}_A\otimes\Omega_B \oplus \Omega_A
\otimes\tilde{\H}_B$. 
As a consequence,  
with the aid of some linear algebra and the fact that $F_A$ and $F_B$
commute with $\Theta$,
it is not difficult to show that
there is a partition in two disjoint sets
$\{1,\ldots,n\}=\alpha_A \cup \alpha_B$ 
along with a suitable choice of the generating
fields such that, for every $f \in {\cal S}(\RR^4)$,
\begin{equation}
\varphi_i(f)\Omega \in \tilde{\H}_A\otimes\Omega_B
\;  {\rm for} \; i\in \alpha_A , \quad
\varphi_i(f)\Omega \in \Omega_A \otimes\tilde{\H}_B 
\; {\rm for} \; i \in\alpha_B .
\label{AB}
\end{equation}

\medskip

Because of equations \ref{AB}, for every $f \in {\cal S}(\RR^4)$
and $i \in \alpha_A$  
one can define a vector $T_i(f) \in \tilde\H_A$
by
\begin{equation}
\varphi_i(f)(\Omega_A \otimes \Omega_B)=: T_i(f) \otimes \Omega_B .
\end{equation}

It follows that
if ${\rm supp}(f) \subset \O$, $f$ real, 
and $X_A \in \F_A(\O')$, $X_B \in \F_B(\O')$,
we get that
\begin{align}
\nonumber
\varphi_i(f)(X_A \Omega_A \otimes X_B\Omega_B) 
& = \varphi_i(f)(X_A \otimes X_B)(\Omega_A \otimes \Omega_B)\\
\nonumber
& = (X_A \otimes X_B)\varphi_i(f)(\Omega_A \otimes \Omega_B)\\
& = X_A T_i(f)\otimes X_B\Omega_B, \; i\in \alpha_A .
\end{align}

By a continuity argument 
(we are assuming $\varphi_i(f)$ to be closed), 
$$\varphi_i(f)(X_A\Omega_A \otimes \xi)=X_A T_i(f)\otimes \xi
\quad
\forall \; \xi \in \H_B.$$

Therefore, for every $T \in {\rm B}(\H_B)$, 
$(I\otimes T)(X_A\Omega_A \otimes X_B \Omega_B)$  
belongs to the domain of $\varphi_i(f)$ and 
\begin{equation}
(I \otimes T) \varphi_i(f) (X_A\Omega_A \otimes X_B \Omega_B)
= \varphi_i(f) (I \otimes T) (X_A\Omega_A \otimes X_B \Omega_B).
\end{equation}
Hence again by continuity we find that, for every $X\in\F(\O ')$,
\begin{equation}
(I \otimes T) \varphi_i(f) X\Omega
= \varphi_i(f) (I \otimes T) X\Omega ,\; i \in \alpha_A .
\label{commutA}
\end{equation}

Similarly, for each $T \in {\rm B}(\H_A)$,
\begin{equation}
(T \otimes I) \varphi_i(f) X\Omega
= \varphi_i(f) (T \otimes I) X\Omega ,\; i \in \alpha_B .
\label{commutB}
\end{equation}

Our next goal is to show that $\F(\O')\Omega$ is a core for $\varphi_i(f)$
for any real $f$ as above and $i =1,\ldots,n$.
This will entail that
$e^{i\varphi_i(f)}\in (I\otimes{\rm B}(\H_B))'= {\rm B}(\H_A)\otimes I$
for every real-valued test function $f$ with compact support
(by arbitrariness of $\O$ in the argument above) and $i \in \alpha_A$,
and similarly 
$e^{i\varphi_i(f)}\in I\otimes{\rm B}(\H_B)$ for $i \in \alpha_B$,
from which it is easy to see that
$\vee_{i \in \alpha_A} \F_i(\O) = \F_A(\O) \otimes I$ and
$\vee_{i \in \alpha_B} \F_i(\O) = I \otimes \F_B(\O)$, $\O\in\K$.

\begin{proposition}
For any $f \in {\cal S}(\RR^4)$ real, $\O\in\K$ and $i=1,\ldots,n$,
$\F(\O)\Omega$ contains a core for $\varphi_i(f)$.
In particular if ${\rm supp}(f)\subset \O$ then $\F(\O')\Omega$ is a core
for $\varphi_i(f)$.
\end{proposition}

\noindent{\it Proof.}
We use some techniques concerning energy-bounds,
cf. \cite[Section 13.1.3]{BaWo}.
Let $N$ be the total number operator acting on $\H=\Gamma(\k)$.
Then $N$ is the closure of $\sum_i N_i$ with $N_i$ the number operator on
$\Gamma(\k_i)$.
Using well known enstimates about free fields
(see \cite[Section X.7]{RS2})
for every real $f$ and $\psi$ in the domain of $N$ we have
\begin{equation}
\|\varphi_i(f)\psi \|
\leq c(f) \|\sqrt{N+I}\psi\|
\leq c(f) \|(N+I)\psi\|
\end{equation}
for some constant $c(f)$ depending only on $f$.
Moreover $\varphi_i(f)$ is essentially self-adjoint on any core for $N$.

We define a self-adjoint operator $H$ as (the closure of) the sum of the
$H_i$, where $H_i$ on $\Gamma(\k_i)$ is the conformal Hamiltonian 
if $\varphi_i(x)$ has vanishing mass and the generator of time translations
otherwise. 
Note that $N^2_i \leq c^2_i H^2_i$, where $c_i$ is the inverse of the mass 
corresponding to $\varphi_i(x)$ if that is different from 0, and equal to 1
otherwise.

It follows that, for $\psi$ in the domain of $H$,
\begin{equation}
\|\varphi_i(f)\psi \|
\leq b(f) \|(H+I)\psi\|
\end{equation}
for some constant $b(f)$.

Thus, since $N$ is essentially self-adjoint on the domain of $H$,
$\varphi_i(f)$ is essentially self-adjoint on any core for $H$.

To complete the proof we only need to show that, for each $\O\in\K$,
$\F(\O)\Omega$ contains a core for $H$.
But this follows from \cite[Appendix]{seba}, 
after noticing that given $\O_1 \subset\subset \O$ then 
$e^{itH}\F(\O_1)e^{-itH} \subset \F(\O)$
for $|t|$ small enough.
\qed

\medskip
Summing up, we have thus proved the following result.

\begin{theorem}
Let $\F$ be the net generated by a finite family of free Hermitian scalar
fields and let $\F=\F_A \otimes \F_B$ be a tensor product decomposition,
then, for a suitable choice $\varphi_1(x),\ldots,$ $\varphi_n(x)$ of the
generating fields for $\F$ and a $k\in \{1,\ldots,n\}$,
$\F_A\otimes I$ is generated by $\varphi_1(x),\ldots$,$\varphi_k(x)$ and
$I\otimes \F_B$ 
by $\varphi_{k+1}(x),\ldots$,$\varphi_n(x)$.
\label{freedec}
\end{theorem}

\medskip
\noindent{\bf Acknowledgements.}
We thank D. R. Davidson, S. Doplicher, J. E. Roberts 
for some useful comments and discussions
at different stages of this research.
Part of this work has been done while R. C. was visiting the 
Department of Mathematics at the University of Oslo.
He thanks the members of the operator algebras group in Oslo for their
warm hospitality
and the EU TMR network ``Non-commutative geometry'' for financial support.


\medskip

\begin{thebibliography}{99}

\bibitem{Ar} Araki H., Symmetries in the theory of local observables
and the choice of the net of local algebras,
{\it Rev. Math. Phys.\/}, {\bf Special Issue} (1992), 1--14.

\bibitem{Barg} Bargmann V.: 
On unitary ray representations of continuous groups, 
{\it Ann. Math.\/} {\bf 59} (1954), 1--46.

\bibitem{BaWo} Baumg\"artel H., Wollenberg M.: 
{\it Causal Nets of Operator Algebras. Mathematical Aspects of Algebraic
Quantum Field Theory}, Akademie Verlag, Berlin 1992. 

\bibitem{BCL} Bertozzini P., Conti R., Longo R.:
Covariant sectors with infinite dimension and positivity of the energy,
{\it Comm. Math. Phys.} {\bf 193} (1998), 471--492.

\bibitem{Borch99} Borchers H.J.: On the revolutionalization of quantum 
field theory by Tomita's modular theory, 
preprint, University of G\"ottingen (1999).

\bibitem{BuDA} Buchholz D., D'Antoni C.:
Phase space properties of charged fields in theories of local observables,
{\it Rev. Math. Phys.} {\bf 7} (1995), 527--557.

\bibitem{BDLR} Buchholz D., Doplicher S., Longo R., Roberts J. E.: A new 
look at Goldstone theorem, {\it Rev. Math.Phys.}, {\bf Special Issue} 
(1992), 47--82.

\bibitem{ext} Buchholz D., Doplicher S., Longo R., Roberts J. E.:
Extension of automorphisms and gauge symmetries, {\it Comm. Math. Phys.} 
{\bf 155} (1993), 123--134.

\bibitem{seba2} Carpi S.: Absence of subsystems for the Haag-Kastler net
generated by the energy-momentum in two dimensional conformal field
theory, {\it Lett. Math. Phys.} {\bf 45} (1998), 259--267.

\bibitem{seba} Carpi S.: Quantum Noether's theorem and conformal field
theory: a study of some models, {\it Rev. Math. Phys.} {\bf 11} (1999),
519--532.

\bibitem{Car} Carpi S.: Classification of subsystems for the
Haag-Kastler nets generated by $c=1$ chiral current algebras,
{\it Lett. Math. Phys.} {\bf 47} (1999), 353--364.

\bibitem{roberto1} Conti R.: On the intrinsic definition of local 
observables, {\it Lett. Math. Phys.} {\bf 35} (1995), 237--250.

\bibitem{roberto} Conti R.: Teoria algebrica dei campi e inclusioni di
algebre di von Neumann, Tesi di Dottorato, 
University of Rome ``Tor Vergata'' (1996).

\bibitem{CoDA} Conti R., D'Antoni C.: Extension of anti--automorphisms
and PCT--symmetry, {\it Rev. Math. Phys.}, to appear.

\bibitem{CDR} Conti R., Doplicher S., Roberts J. E.:
Superselection theory for subsystems, preprint, 
University of Rome ``La Sapienza'' (1999).

\bibitem{DA} D'Antoni C.: Technical properties of the quasi-local algebra, 
in \cite{Kast}.

\bibitem{DALR} D'Antoni C., Longo R., Radulescu F.:
Conformal nets, maximal temperature and models from free probability,
{\it J. Operator Th.}, to appear.

\bibitem{Dav} Davidson D. R.: 
Classification of subsystems of local algebras, 
{\it Ph.D. Thesis}, University of California at Berkeley (1993).

\bibitem{Dix1} Dixmier J.: {\it Les $C^*$-alg\`ebres et leurs
repr\'esentations,} Gauthier Villars, Paris 1964.

\bibitem{Dix2} Dixmier J.: {\it Les alg\`ebres d'op\'erateurs dans
l'espace hilbertien,} Gauthier Villars, Paris 1969.

\bibitem{Do82} Doplicher S.: Local aspects of superselection rules, 
{\it Comm. Math. Phys.} {\bf 85} (1982), 73--85. 

\bibitem{Dop} Doplicher S.:
Progress and problems in algebraic quantum
field theory, in S. Albeverio et al. (eds.), Ideas and Methods in Quantum
and Statistical Physics Vol. 2.

\bibitem{DHR4} Doplicher S., Haag R., Roberts J. E.: LocaI observable and 
particle statistics II, {\it Comm. Math. Phys.} {\bf 35} (1974), 49--85.

\bibitem{DoLo84} Doplicher S., Longo R.: Standard and split inclusions 
of von Neumann algebras, {\it Invent. Math.} {\bf 75} (1984), 493--536. 

\bibitem{DP} Doplicher S., Piacitelli G.: In preparation.

\bibitem{DoRo90} Doplicher S., Roberts J.E.: Why there is a field
algebra with a compact gauge group describing the superselection structure
structure in particle physics, {\it Comm. Math. Phys.\/}
{\bf 131} (1990), 51--107.

\bibitem{Dr} Driessler W.: Duality and absence of locally generated sectors
for CCR-type algebras, {\it Comm. Math. Phys.\/} {\bf 70} (1979),
213--220.

\bibitem{GeKa} Ge L., Kadison R. V.: On tensor products of von Neumann 
algebras, {\it Invent. Math.} {\bf 123} (1996), 453--466.

\bibitem{GuLo1} Guido D., Longo R.:
An algebraic spin and statistics theorem. I,
{\it Comm. Math. Phys.\/} {\bf 172} (1995), 517--533.

\bibitem{GLRV} Guido D., Longo R., Roberts J. E., Verch R.:
Charged sectors, spin and statistics in quantum field theory on curved
spacetime, preprint, University of Rome ``Tor Vergata'' (1999).

\bibitem{Ha} Haag R.: {\it Local Quantum Physics}, 2nd ed.\,  
Springer-Verlag, New York Berlin Heidelberg 1996.

\bibitem{KR2} Kadison R.V., Ringrose J.R.: {\it Fundamentals of the theory 
of operator algebras.} Volume II. Academic Press, New York, 1986.

\bibitem{Kast} Kastler D. ed.:{\it The algebraic theory of superselection 
sectors.} World Scientific, Singapore, 1990. 

\bibitem{KLM} Kawahigashi Y., Longo R., M\"{u}ger M.: Multi-interval 
subfactor and modularity of representations in conformal field theory,
preprint, University of Rome ``Tor Vergata'' (1999).

\bibitem{Kuc} Kuckert B.:
A new approach to Spin and Statistics,
{\it Lett. Math. Phys.} {\bf 35} (1995), 319--331.

\bibitem{Kuc2} Kuckert B.: Spin \& Statistics, localization regions,
and modular symmetries in quantum field theory, 
{\it Ph.D. Thesis}, Universit\"{a}t Hamburg (1998).

\bibitem{LaSc1} Langerholc J., Schroer B.: 
On the structure of the von Neumann algebras generated by local functions
of the free Bose field,
{\it Comm. Math. Phys.\/} {\bf 1} (1965), 215--239.

\bibitem{LaSc2} Langerholc J., Schroer B.: 
Can current operators determine a complete theory?,  
{\it Comm. Math. Phys.\/} {\bf 4} (1967), 123--136.

\bibitem{longo} Longo R.: Private communication.

\bibitem{LR} Longo R., Rehren K.-H.: Nets of subfactors, 
{\it Rev. Math. Phys.} {\bf 7} (1995), 567--597.

\bibitem{Mu2} M\"uger M.:
Quantum double actions on operator algebras and orbifold quantum field
theories,
{\it Comm. Math. Phys.} {\bf 191} (1998), 137--181.

\bibitem{Mu} M\"uger M.:
On charged fields with group symmetry and degeneracies of Verlinde's
matrix $S$, {\it Ann. Inst. H. Poincar\'e}
{\bf 71} (1999), 359--394.

\bibitem{RS2} Reed M., Simon B.: {\it Methods of modern mathematical
physics} vol. II, Academic Press,  New York 1975.

\bibitem{StZs} Str\u{a}til\u{a} S., Zsido L.:
The commutation theorem for tensor products over von Neumann subalgebras,
{\it J. Funct. Anal.} {\bf 165} (1999), 293--346.

\bibitem{TW} Thomas L.J. III, Wichmann E. H.: Standard forms of local nets 
in quantum field theory, {\it J. Math. Phys.} {\bf 39} (1998), 2643--2681.

\bibitem{Wic} Wichmann E. H.: On systems of local operators and the duality 
condition, {\it J. Math. Phys.} {\bf 24} (1983), 1633--1644.

\bibitem{Wig} Wigner E.: 
On unitary representations of the inhomogeneous Lorentz group,
{\it Ann. Math.\/} {\bf 40} (1939), 149--204.

\end{thebibliography}
\end{document}